\documentclass[12pt]{article}
\usepackage{amsmath}
\usepackage{amsfonts}
\usepackage{amssymb}
\usepackage{longtable}

\newtheorem{theorem}{Theorem}{\bf}{\it}
\newtheorem{lemma}{Lemma}{\bf}{\it}
\newtheorem{proposition}{Proposition}{\bf}{\it}

\newcommand{\id}{\mathop\text{\rm id}\nolimits}
\newcommand{\pr}{\mathop\text{\rm pr}\nolimits}
\newcommand{\Ric}{\mathop{{\rm Ric}}\nolimits}
\newcommand{\tRic}{\mathop{\widetilde{\rm Ric}}\nolimits}

\begin{document}

\title{About the classification of the holonomy algebras of Lorentzian manifolds}
\author{Anton S. Galaev}

%\author{Anton S. Galaev}
%\address{Department of Mathematics and Statistics, Faculty of Science, Masaryk University in Brno,
%Kotl\'a\v rsk\'a~2, 611~37 Brno,
%Czech Republic\\ \ead{galaev@math.muni.cz}}

\maketitle

\begin{abstract}
The classification of the holonomy algebras of Lorentzian
manifolds can be reduced to the classification of irreducible
subalgebras $\mathfrak{h}\subset\mathfrak{so}(n)$ that are spanned
by  the images of linear maps from $\mathbb{R}^n$ to
$\mathfrak{h}$ satisfying an identity similar to the Bianchi one.
T. Leistner found all such subalgebras and it turned out that the
obtained list coincides with the list of irreducible holonomy
algebras of Riemannian manifolds. The natural problem is to give a
simple direct proof to this fact.  We give such proof for the case
of semisimple not simple Lie algebras $\mathfrak{h}$.

\end{abstract}

\section{Introduction}

M.~Berger \cite{Ber,Besse} classified possible connected
irreducible holonomy groups $H\subset {\rm SO}(n)$ of not locally
symmetric Riemannian manifolds using the representation theory. It
turned out that these groups  act transitively on the unite sphere
of the tangent space. J.~Simens \cite{Simens} and recently in a
simple geometric way C.~Olmos \cite{Olmos} proved this result
directly.

The classification of the holonomy algebras (i.e. the Lie algebras
of the holonomy groups) of Lorentzian manifolds can be reduced to
the classification of irreducible weak-Berger subalgebras
$\mathfrak{h}\subset\mathfrak{so}(n)$, i.e. subalgebras
$\mathfrak{h}\subset\mathfrak{so}(n)$ that are spanned by the
images of linear maps from the space
$$\mathcal{P}(\mathfrak{h})=\{P\in(\mathbb{R}^n)^*\otimes
\mathfrak{h}|(P(X)Y,Z)+(P(Y)Z,X)+(P(Z)X,Y)=0,\,\,
X,Y,Z\in\mathbb{R}^n\}.$$ It is easy to see that if
$\mathfrak{h}\subset\mathfrak{so}(n)$ is the  holonomy algebras of
a Riemannian manifold, then it is a weak-Berger algebra. The
inverse statement is absolutely not obvious, nevertheless it is
true and it is proven by Th.~Leistner in \cite{Leistner}.

If $n$ is even and $\mathfrak{h}\subset\mathfrak{so}(n)$ is of
complex type, i.e. $\mathfrak{h}\subset\mathfrak{u}(\frac{n}{2})$,
then it can be shown that $\mathcal{P}(\mathfrak{h})\simeq
(\mathfrak{h}\otimes\mathbb{C})^{(1)}$, where
$(\mathfrak{h}\otimes\mathbb{C})^{(1)}$ is the first prolongation
of the subalgebra
$\mathfrak{h}\otimes\mathbb{C}\subset\mathfrak{gl}(\frac{n}{2},\mathbb{C})$
(cf. \cite{Leistner} and \cite{onecomp}). Using that and the
classification of irreducible representations with non-trivial
prolongation, Leistner showed that if
$\mathfrak{h}\subset\mathfrak{u}(\frac{n}{2})$ is a weak-Berger
subalgebra, then it is the holonomy algebra of a Riemannian
manifold.

The situation when $\mathfrak{h}\subset\mathfrak{so}(n)$ is of
real type (i.e. not of complex type) is much more difficult. In
this case Leistner considered the complexification
$\mathfrak{h}\otimes\mathbb{C}\subset\mathfrak{so}(n,\mathbb{C})$,
which is irreducible. He used the classification of irreducible
representations of complex semisimple Lie algebras, found a
criteria in terms of weights for such representation
$\mathfrak{h}\otimes\mathbb{C}\subset\mathfrak{so}(n,\mathbb{C})$
to be a weak-Berger algebra and considered case by case simple Lie
algebras $\mathfrak{h}\otimes\mathbb{C}$, and then semisimple Lie
algebras (the problem is reduced to the semisimple Lie algebras of
the form $\mathfrak{sl}(2,\mathbb{C})\oplus\mathfrak{k}$, where
$\mathfrak{k}$ is simple, and again different possibilities for
$\mathfrak{k}$ were considered).

We consider the case of semisimple not simple irreducible
subalgebras $\mathfrak{h}\subset\mathfrak{so}(n)$ with irreducible
complexification
$\mathfrak{h}\otimes\mathbb{C}\subset\mathfrak{so}(n,\mathbb{C})$.
In a simple way we show that it is enough to treat the case when
$\mathfrak{h}\otimes\mathbb{C}=\mathfrak{sl}(2,\mathbb{C})\oplus\mathfrak{k}$,
where $\mathfrak{k}\subsetneq\mathfrak{sp}(2m,\mathbb{C})$ is a
proper irreducible subalgebra, and the representation space is the
tensor product $\mathbb{C}^2\otimes\mathbb{C}^{2m}$. We show that
in this case $\mathcal{P}(\mathfrak{h})$ coincides with
$\mathbb{C}^2\otimes\mathfrak{g}_1$, where $\mathfrak{g}_1$ is the
first Tanaka prolongation the non-positively graded Lie algebra
$$\mathfrak{g}_{-2}\oplus\mathfrak{g}_{-1}\oplus\mathfrak{g}_0,$$ where $\mathfrak{g}_{-2}=\mathbb{C},\quad
\mathfrak{g}_{-1}=\mathbb{C}^{2m},\quad
\mathfrak{g}_0=\mathfrak{k}\oplus\mathbb{C}
\id_{\mathbb{C}^{2m}}$, and the grading is defined by the element
$-\id_{\mathbb{C}^{2m}}$. We prove that if
$\mathcal{P}(\mathfrak{h})$ is non-zero, then $\mathfrak{g}_1$ is
isomorphic to $\mathbb{C}^{2m}$, the second Tanaka prolongation
$\mathfrak{g}_2$ is isomorphic to $\mathbb{C}$, and
$\mathfrak{g}_3=0$. Then, the full Tanaka prolongation defines the
simple $|2|$-graded complex Lie algebra
$$\mathfrak{g}_{-2}\oplus\mathfrak{g}_{-1}\oplus\mathfrak{g}_0\oplus\mathfrak{g}_{1}\oplus\mathfrak{g}_2.$$
It is well known that such Lie algebra defines (up the duality) a
simply connected Riemannian symmetric space; the holonomy algebra
of this space coincides with
$\mathfrak{h}\subset\mathfrak{so}(n)$. Thus, if the subalgebra
$\mathfrak{h}\subset\mathfrak{so}(n)$ is semisimple and not
simple, and $\mathcal{P}(\mathfrak{h})\neq 0$, then we indicate a
Riemannian manifold with the holonomy algebra
$\mathfrak{h}\subset\mathfrak{so}(n)$.

More details about the holonomy algebras of Lorentzian manifolds
can be found in \cite{Gal5,ESI}.

\section{Holonomy algebras of Riemannian manifolds}\label{secRiemannian}

Irreducible holonomy algebras
$\mathfrak{h}\subset\mathfrak{so}(n)$ of not locally symmetric
Riemannian manifolds are exhausted by $\mathfrak{so}(n)$,
$\mathfrak{u}(\frac{n}{2})$, $\mathfrak{su}(\frac{n}{2})$,
$\mathfrak{sp}(\frac{n}{4})\oplus\mathfrak{sp}(1)$,
$\mathfrak{sp}(\frac{n}{4})$, $G_2\subset\mathfrak{so}(7)$ and
$\mathfrak{spin}(7)\subset\mathfrak{so}(8)$. This list (up to some
corrections) obtained M. Berger \cite{Ber,Besse}. Berger
classified irreducible subalgebras
$\mathfrak{h}\subset\mathfrak{so}(n)$ spanned by the images of the
maps from the space
 $$\mathcal{R}(\mathfrak{h})=\{R\in \Lambda^2
(\mathbb{R}^n)^*\otimes\mathfrak{h}|R(X,Y)Z+R(Y, Z)X+R(Z, X)Y=0
\text{ for all } X,Y,Z\in \mathbb{R}^n\}$$  of algebraic curvature
tensors of type $\mathfrak{h}$ under the condition that the space
$$\mathcal{R}^\nabla(\mathfrak{h})=\{S\in
(\mathbb{R}^n)^*\otimes\mathcal{R}(\mathfrak{h})|S_X(Y,Z)+S_Y(Z,X)+S_Z(X,Y)=0
\text{ for all } X,Y,Z\in \mathbb{R}^n\}$$ of algebraic covariant
derivatives of the curvature tensors of type $\mathfrak{h}$ is not
trivial. Berger used the classification of irreducible
representations of compact Lie groups. The connected Lie subgroups
of ${\rm SO}(n)$ corresponding to the above subalgebras of
$\mathfrak{so}(n)$ mostly exhaust
  groups of isometries acting transitively on the unite sphere of dimension $n-1$, and the result of Berger can be reformulated in
the following form: if the irreducible holonomy group of a
Riemannian manifold $(M,g)$ does not act transitively on the unite
sphere of the tangent space, then $(M,g)$ is locally symmetric. A
direct proof of this statement obtained in algebraic way J.~Simens
\cite{Simens}, and recently an elegant geometric proof obtained
C.~Olmos \cite{Olmos}.

The spaces $\mathcal{R}(\mathfrak{h})$ for the irreducible
holonomy algebras of Riemannian manifolds
$\mathfrak{h}\subset\mathfrak{so}(n)$ are computed by
D.~V.~Alekseevsky in \cite{Al}.  For $R\in
\mathcal{R}(\mathfrak{h})$ define its Ricci tensor by
$$\Ric(R)(X,Y)={\rm tr}(Z\mapsto R(Z,X)Y),$$ $X,Y\in \mathbb{R}^n$.
 The space  $\mathcal{R}(\mathfrak{h})$ admits the following decomposition into
$\mathfrak{h}$-modules:
$$\mathcal{R}(\mathfrak{h})=\mathcal{R}_0(\mathfrak{h})\oplus\mathcal{R}_1(\mathfrak{h})\oplus\mathcal{R}'(\mathfrak{h}),$$ where
$\mathcal{R}_0(\mathfrak{h})$ consists of the curvature tensors
with zero Ricci curvature, $\mathcal{R}_1(\mathfrak{h})$ consists
of tensors annihilated by $\mathfrak{h}$ (this space is either
zero or one-dimensional), $\mathcal{R}'(\mathfrak{h})$ is the
complement to these two spaces. If
$\mathcal{R}(\mathfrak{h})=\mathcal{R}_1(\mathfrak{h})$, then any
Riemannian manifold with the holonomy algebra
$\mathfrak{h}\subset\mathfrak{so}(n)$ is locally symmetric. Such
subalgebra $\mathfrak{h}\subset\mathfrak{so}(n)$ is called {\it a
symmetric Berger algebra}.  The holonomy algebras of irreducible
Riemannian symmetric spaces are exhausted by $\mathfrak{so}(n)$,
$\mathfrak{u}(\frac{n}{2})$,
$\mathfrak{sp}(\frac{n}{4})\oplus\mathfrak{sp}(1)$ and by
symmetric Berger algebras $\mathfrak{h}\subset\mathfrak{so}(n)$.

It is known that simply connected indecomposable symmetric
Riemannian manifolds $(M,g)$ are in one-two-one correspondence
with simple $\mathbb{Z}_2$-graded Lie algebras
$\mathfrak{g}=\mathfrak{h}\oplus\mathbb{R}^n$ such that
$\mathfrak{h}\subset\mathfrak{so}(n)$. The subalgebra
$\mathfrak{h}\subset\mathfrak{so}(n)$ coincides with the holonomy
algebra of $(M,g)$. The space $(M,g)$ can be reconstructed using
its holonomy algebra $\mathfrak{h}\subset\mathfrak{so}(n)$ and the
value $R\in\mathcal{R}(\mathfrak{h})$ of the curvature tensor of
$(M,g)$ at a point. For that define the structure of the Lie
algebra on the vector space
$\mathfrak{g}=\mathfrak{h}\oplus\mathbb{R}^n$ in the following
way:
$$[A,B]=[A,B]_{\mathfrak{h}},\quad [A,X]=AX,\quad [X,Y]=R(X,Y),\quad
A,B\in\mathfrak{h},\,X,Y\in\mathbb{R}^n.$$ Then $M=G/H$, where $G$
is the simply connected Lie group corresponding to the Lie algebra
$\mathfrak{g}$, and $H\subset G$ is the connected Lie subgroup
corresponding to the subalgebra $\mathfrak{h}\subset\mathfrak{g}$.

If the symmetric space is quaternionic-K\"ahlerian, then
$\mathfrak{h}=\mathfrak{sp}(1)\oplus\mathfrak{f}\subset\mathfrak{so}(4k)$,
where $n=4k$ and $\mathfrak{f}\subset\mathfrak{sp}(k)$. The
complexification of $\mathfrak{h}\oplus\mathbb{R}^{4k}$ is equal
to
$(\mathfrak{sl}(2,\mathbb{C})\oplus\mathfrak{k})\oplus(\mathbb{C}^2\otimes
\mathbb{C}^{2k})$, where
$\mathfrak{k}=\mathfrak{f}\otimes\mathbb{C}\subset\mathfrak{sp}(2k,\mathbb{C})$.
Let $e_1,e_2$ be the standard basis of $\mathbb{C}^2$, and  let
$$F=\left(\begin{matrix}0&0\\1&0\end{matrix}\right),\quad H=\left(\begin{matrix}1&0\\0&-1\end{matrix}\right),\quad
E=\left(\begin{matrix}0&1\\0&0\end{matrix}\right)$$ be the  basis
of $\mathfrak{sl}(2,\mathbb{C})$. We obtain the following
$\mathbb{Z}$-grading of $\mathfrak{g}\otimes\mathbb{C}$:
$$\mathfrak{g}\otimes\mathbb{C}=\mathfrak{g}_{-2}\oplus\mathfrak{g}_{-1}\oplus\mathfrak{g}_0\oplus\mathfrak{g}_1\oplus\mathfrak{g}_2=\mathbb{C}
F\oplus e_2\otimes
\mathbb{C}^{2k}\oplus(\mathfrak{k}\oplus\mathbb{C} H)\oplus
e_1\otimes\mathbb{C}^{2k}\oplus\mathbb{C} E.$$ Conversely, any
such simple $\mathbb{Z}$-graded Lie algebra defines up to the
duality a simply connected quaternionic-K\"ahlerian symmetric
space.

\section{Weak curvature tensors}

The spaces $\mathcal{P}(\mathfrak{h})$ are computed in
\cite{onecomp}. Let $\mathfrak{h}\subset\mathfrak{so}(n)$ be an
irreducible subalgebra. There exists the decomposition
$$\mathcal{P}(\mathfrak{h})=\mathcal{P}_0(\mathfrak{h})\oplus\mathcal{P}_1(\mathfrak{h}),$$ where $\mathcal{P}_0(\mathfrak{h})$ is the kernel
of the $\mathfrak{h}$-equivariant map
$$\tRic:\mathcal{P}(\mathfrak{h})\to\mathbb{R}^n, \qquad \tRic(P)=\sum_{i=1}^nP(e_i)e_i$$
($e_1,...,e_n$ is an orthogonal basis of $\mathbb{R}^n$), and
$\mathcal{P}_1(\mathfrak{h})$ is the orthogonal complement of
$\mathcal{P}_0(\mathfrak{h})$ in $\mathcal{P}(\mathfrak{h})$. The
space $\mathcal{P}_1(\mathfrak{h})$ is either trivial or it is
isomorphic to $\mathbb{R}^n$. If $n\geq 4$, then
$\mathcal{P}_0(\mathfrak{h})\neq 0$ if and only if
$\mathcal{R}_0(\mathfrak{h})\neq 0$. Next,
$\mathcal{P}_1(\mathfrak{h})\simeq\mathbb{R}^n$  if and only if
$\mathcal{R}_1(\mathfrak{h})\simeq \mathbb{R}$. For the symmetric
Berger algebras it holds
$\mathcal{P}_1(\mathfrak{h})\simeq\mathbb{R}^n$ and
$\mathcal{P}_0(\mathfrak{h})=0$.

\section{Tanaka prolongations}

Consider a $\mathbb{Z}$-graded Lie algebra of the form
$$\mathfrak{g}_{-2}\oplus\mathfrak{g}_{-1}\oplus\mathfrak{g}_0.$$

For $k\geq 1$, the  $k$-th Tanaka prolongation is defined by the
induction
$$\mathfrak{g}_k=\{u\in (\mathfrak{g}^*_{-2}\otimes\mathfrak{g}_{k-2})\oplus
(\mathfrak{g}^*_{-1}\otimes\mathfrak{g}_{k-1})|u([X,Y])=[u(X),Y]+[X,u(Y)],\,\,X,Y\in
\mathfrak{g}_{-2}\oplus\mathfrak{g}_{-1}\}.$$ Let $k\geq 1$, and
$l\geq 0$. For $u\in\mathfrak{g}_k$ and $v\in\mathfrak{g}_l$
define Lie brackets $[u,v]\in\mathfrak{g}_{k+l}$ by the condition
$$[u,v]X=[[u,X],v]+[u,[v,X]],\quad X\in \mathfrak{g}_{-2}\oplus\mathfrak{g}_{-1};$$
the Lie brackets of $u\in\mathfrak{g}_k$ and $X\in
\mathfrak{g}_{-2}\oplus\mathfrak{g}_{-1}$ are defined as
$[u,X]=Xu$. This gives the structure of a Lie algebra on the
vector space $\oplus_{k=-2}^\infty\mathfrak{g}_k$.

Let $\mathfrak{k}\subset\mathfrak{sp}(2m,\mathbb{C})$ be a
subalgebra, $m\geq 2$. Consider the Lie algebra
$$\mathfrak{g}_{-2}\oplus\mathfrak{g}_{-1}\oplus\mathfrak{g}_0,\quad
\mathfrak{g}_{-2}=\mathbb{C} F,\quad
\mathfrak{g}_{-1}=\mathbb{C}^{2m},\quad\mathfrak{g}_0=\mathfrak{k}\oplus\mathbb{C}
H$$ with the non-zero Lie brackets
$$[X,Y]=\Omega(X,Y)F,\quad [A,X]=AX, \quad [A,B]=[A,B]_{\mathfrak{k}},\quad [H,X]=-X,\quad
[H,F]=-2F,$$ where $X,Y\in\mathbb{C}^{2m}$, $A,B\in\mathfrak{k},$
and $\Omega$ is the symplectic form on $\mathbb{C}^{2m}$.

\begin{lemma}\label{lemTan}   It holds
$$\mathfrak{g}_1=\{\varphi\in\mathfrak{g}_{-1}^*\otimes\mathfrak{g}_0|\exists A\in \mathfrak{g}_{-1},\,
\varphi(X)Y-\varphi(Y)X=\Omega(X,Y)A,\,
X,Y\in\mathfrak{g}_{-1}\}.$$  If
$\mathfrak{k}\subsetneq\mathfrak{sp}(2m,\mathbb{C})$ is a proper
irreducible subalgebra and $\mathfrak{g}_1\neq 0$, then
$\mathfrak{g}_1\simeq\mathbb{C}^{2m}$, $\mathfrak{g}_2\simeq
\mathbb{C}$, and $\mathfrak{g}_3=0$. The Lie algebra
$$\mathfrak{g}_{-2}\oplus\mathfrak{g}_{-1}\oplus\mathfrak{g}_0\oplus\mathfrak{g}_1\oplus\mathfrak{g}_2$$ is simple.

\end{lemma}

{\bf Proof.}  Let $u=\psi+\varphi,$ where $\psi\in
\mathfrak{g}^*_{-2}\otimes\mathfrak{g}_{-1}$, and
$\varphi\in\mathfrak{g}^*_{-1}\otimes\mathfrak{g}_{0}$. The
condition $u\in\mathfrak{g}_1$ is equivalent to the equations
$$[\varphi(X),F]=\Omega(\psi(F),Y)F,\quad
\varphi(X)Y-\varphi(Y)X=\Omega(X,Y)\psi(F).$$ The first statement
of the lemma is that the second equation implies the first one.

Let us denote $\mathbb{C}^{2m}$ by $V$. First suppose that
$\mathfrak{k}=\mathfrak{sp}(V)$. Let us find $\mathfrak{g}_1$. We
have the following isomorphisms of the $\mathfrak{sp}(V)$-modules:
$\mathfrak{g}^*_{-2}\otimes\mathfrak{g}_{-1}\simeq V$, and
$$\mathfrak{g}^*_{-1}\otimes\mathfrak{g}_{0}\simeq V\otimes (\mathfrak{sp}(V)\oplus
\mathbb{C})=V\oplus(V\oplus V_{3\pi_1}\oplus V_{\pi_1+\pi_2}),$$
where $V_\Lambda$ denotes the irreducible
$\mathfrak{sp}(V)$-module with the highest weight $\Lambda$. By
the definition, the intersection of $\mathfrak{g}_1$ and
$\mathfrak{g}^*_{-1}\otimes\mathfrak{g}_{0}$ coincides with
$$(\mathfrak{sp}(V)\oplus\mathbb{C} H)^{(1)}=(\mathfrak{sp}(V))^{(1)}=\odot^3V\simeq
V_{3\pi_1}.$$ Clearly, the intersection of $\mathfrak{g}_1$ and
$\mathfrak{g}^*_{-2}\otimes\mathfrak{g}_{-1}$ is trivial.
Consequently, if $\mathfrak{g}_1$ is different from
$\mathfrak{sp}(V)^{(1)}$, then $\mathfrak{g}_1$ contains a
submodule isomorphic to $V$. Any $\mathfrak{sp}(V)$-equivariant
map from $V$ to
$(\mathfrak{g}^*_{-2}\otimes\mathfrak{g}_{-1})\oplus
(\mathfrak{g}^*_{-1}\otimes\mathfrak{g}_{0})$ is of the form
$$Z\mapsto\psi^Z+\varphi^Z,\quad \psi^Z(F)=aZ,\quad
\varphi^Z(Y)=b\Omega(Z,Y)H+cZ\odot Y,$$ where
$a,b,c\in\mathbb{R}$, and $Z\odot Y\in \mathfrak{sp}(V)$ is
defined as
$$(Z\odot Y)X=\Omega(Z,X)Y+\Omega(Y,X)Z.$$ The second equation on
$\mathfrak{g}_1$ takes the form
$$-b\Omega(Z,X)Y+b\Omega(Z,Y)X+c(\Omega(Y,Z)X-\Omega(X,Z)Y+2\Omega(Y,X)Z)=a\Omega(X,Y)Z.$$
This equation should hold for all $X,Y,Z\in V$, and it is
equivalent to $b=-c=-\frac{1}{2}a$ (since $\dim V\geq 4$). The
second equation on $\mathfrak{g}_1$ takes the form
$$-2b\Omega(Z,Y)=a\Omega(Z,Y)$$ and it follows from the first one.
Thus the orthogonal complement to
$(\mathfrak{sp}(2m,\mathbb{C}))^{(1)}$ in $\mathfrak{g}_1$ is
isomorphic to $V$, and the isomorphism is given by
$$Z\in V\mapsto \psi^Z+\varphi^Z,\quad \psi^Z(F)=2 Z,\quad
\varphi^Z(Y)=-\Omega(Z,Y)H+Z\odot Y,\quad Y\in V.$$

Let $\mathfrak{k}\subsetneq\mathfrak{sp}(V)$ be a proper
irreducible subalgebra. It is clear that
$$\mathfrak{g}_1=((\mathfrak{g}^*_{-2}\otimes\mathfrak{g}_{-1})\oplus
(\mathfrak{g}^*_{-1}\otimes\mathfrak{g}_{0}))\cap(\mathfrak{sp}(V)\oplus\mathbb{C}
H)_1,$$ and
$\mathfrak{h}^{(1)}=(\mathfrak{g}^*_{-1}\otimes\mathfrak{g}_{0})\cap
\mathfrak{sp}(V)^{(1)}$. It is known that $\mathfrak{h}^{(1)}=0$.
Consequently, if $\mathfrak{g}_1\neq 0$, than $\mathfrak{g}_1$ is
isomorphic to $V$ and it is included diagonally into
$V\oplus\mathfrak{sp}(V)^{(1)}$.

Consider the full Tanaka prolongation
$\mathfrak{g}=\oplus_{i=-2}^\infty\mathfrak{g}_i$. Let
$\mathfrak{g}^0=\oplus_{i=0}^\infty\mathfrak{g}_i\subset\mathfrak{g}$.
We claim that $\mathfrak{g}$ is a primitive $\mathbb{Z}$-graded
Lie algebra, i.e. $\mathfrak{g}^0\subset\mathfrak{g}$ is a maximal
graded subalgebra and $\mathfrak{g}^0$ contains no graded ideals
of $\mathfrak{g}$ except $\{0\}$. Indeed, suppose that there
exists a subalgebra $\tilde{\mathfrak{g}}\subset\mathfrak{g}$ such
that $\mathfrak{g}^0\subsetneq\tilde{\mathfrak{g}}$. Then
$aF+X\in\tilde{\mathfrak{g}}$ for some $a\in\mathbb{R}$,
$X\in\mathfrak{g}_{-1}$. If $a\neq 0$, then taking
$u\in\mathfrak{g}_1$, we get $0\neq
u(F)\in\tilde{\mathfrak{g}}\cap\mathfrak{g}_{-1}$, i.e. we may
assume that there exists non-zero $X\in\mathfrak{g}_{-1}$ such
that $X\in\tilde{\mathfrak{g}}$. Since $\mathfrak{g}_0$ acts on
$\mathfrak{g}_{-1}$ irreducible, we get
$\mathfrak{g}_{-1}\subset\tilde{\mathfrak{g}}$. Finally,
$[\mathfrak{g}_{-1},\mathfrak{g}_{-1}]=\mathfrak{g}_{-2}$, i.e.
$\mathfrak{g}_{-2}\subset\tilde{\mathfrak{g}}$ and
$\tilde{\mathfrak{g}}=\mathfrak{g}$. Suppose now that
$\tilde{\mathfrak{g}}=\oplus_{i=0}^\infty\tilde{\mathfrak{g}}_i\subset\mathfrak{g}^0$
is a graded ideal. For $X\in\mathfrak{g}_{-1}$ and
$\xi\in\tilde{\mathfrak{g}}_0$ it holds
$[\xi,X]\in\mathfrak{g}_{-1}$. On the other hand,
$[\xi,X]\in\tilde{\mathfrak{g}}$, and we get $[\xi,X]=0$ for all
$X\in\mathfrak{g}_{-1}$. This implies $\tilde{\mathfrak{g}}_0=0$.
In the same way it can be shown that $\tilde{\mathfrak{g}}_k=0$
for all $k\geq 2$. Thus, $\mathfrak{g}$ is a primitive
$\mathbb{Z}$-graded Lie algebra. If $\mathfrak{g}$ is infinite
dimensional, then from \cite[Th. 6.1]{Gul} it follows that
$\mathfrak{g}_0=\mathfrak{sp}(V)\oplus\mathbb{C} H$, which gives a
contradiction, since we assume that
$\mathfrak{k}\subsetneq\mathfrak{sp}(V)$ is a proper subalgebra.
Thus, $\mathfrak{g}$ is of finite dimension. Since the element
$H\in\mathfrak{g}_0$ defines the $\mathbb{Z}$-grading of
$\mathfrak{g}$, any ideal $\mathfrak{t}\subset\mathfrak{g}$ is
graded. As in the above claim it can be shown that either
$\mathfrak{t}=\mathfrak{g}$ or $\mathfrak{t}=0$, i.e.
$\mathfrak{g}$ is a simple Lie algebra. For the Killing form of a
$\mathbb{Z}$-graded Lie algebra it holds
$b(\mathfrak{g}_k,\mathfrak{g}_{l})=0$ for $k\neq -l$. This shows
that $\mathfrak{g}_2\simeq\mathbb{C}$ and $\mathfrak{g}_3=0$. The
lemma is proved. $\Box$

\section{Semisimple not simple weak-Berger algebras}

\begin{theorem} Let $\mathfrak{h}\subset\mathfrak{so}(n)$ be a semisimple not simple
irreducible subalgebra of real type. If
$\mathcal{P}(\mathfrak{h})\neq 0$, then
$\mathfrak{h}\subset\mathfrak{so}(n)$ is the holonomy algebra of a
symmetric Riemannian space.
\end{theorem}

{\bf Proof.} From the assumption of the theorem it follows that
the complexified representation
$\mathfrak{h}\otimes\mathbb{C}\subset\mathfrak{so}(n,\mathbb{C})$
is irreducible. Since $\mathfrak{h}\otimes\mathbb{C}$ is
semisimple and not simple, it can be decomposed into the direct
sum of two ideals,
$\mathfrak{h}\otimes\mathbb{C}=\mathfrak{h}_1\oplus\mathfrak{h}_2$.
The representation of $\mathfrak{h}_1\oplus\mathfrak{h}_2$ on
$\mathbb{C}^n$ must be of the from of the tensor product,
$\mathbb{C}^n=\mathbb{C}^{n_1}\otimes\mathbb{C}^{n_2}$, where
$\mathfrak{h}_1\subset\mathfrak{gl}(n_1,\mathbb{C})$,
$\mathfrak{h}_2\subset\mathfrak{gl}(n_2,\mathbb{C})$ are
irreducible. Since
$\mathfrak{h}_1\oplus\mathfrak{h}_2\subset\mathfrak{so}(n,\mathbb{C})$,
it holds that either
$\mathfrak{h}_1\subset\mathfrak{so}(n_1,\mathbb{C})$,
$\mathfrak{h}_2\subset\mathfrak{so}(n_2,\mathbb{C})$, $n_1,n_2\geq
3$ or $\mathfrak{h}_1\subset\mathfrak{sp}(n_1,\mathbb{C})$,
$\mathfrak{h}_2\subset\mathfrak{sp}(n_2,\mathbb{C})$, $n_1,n_2\geq
2$. In \cite{onecomp} it is shown in a simple way that
$\mathcal{P}(\mathfrak{so}(n_1,\mathbb{C})\oplus\mathfrak{so}(n_2,\mathbb{C}))\simeq\mathbb{C}^n$,
and if $n_1,n_2\geq 3$, then
$\mathcal{P}(\mathfrak{sp}(n_1,\mathbb{C})\oplus\mathfrak{sp}(n_2,\mathbb{C}))\simeq\mathbb{C}^n$.
This implies that if $\mathfrak{h}_1\oplus\mathfrak{h}_2$ is a
proper irreducible subalgebra of
$\mathfrak{so}(n_1,\mathbb{C})\oplus\mathfrak{so}(n_2,\mathbb{C})$
or of
$\mathfrak{sp}(n_1,\mathbb{C})\oplus\mathfrak{sp}(n_2,\mathbb{C})$
with $n_1,n_2\geq 3$, then
$\mathcal{P}(\mathfrak{h}_1\oplus\mathfrak{h}_2)=0$. Note that the
holonomy algebras of the Riemannian symmetric spaces
 $${\rm
SO}(n_1+n_2)/({\rm SO}(n_1)\times {\rm SO}(n_2)),\, n_1,n_2\geq
3,\quad {\rm Sp}(n_1+n_2)/({\rm Sp}(n_1)\times {\rm Sp}(n_2)),\,
n_1,n_2\geq 1$$ are respectively
$\mathfrak{so}(n_1)\oplus\mathfrak{so}(n_2)$ and
$\mathfrak{sp}(n_1)\oplus\mathfrak{sp}(n_2)$ \cite{Besse}.

Thus, we are left with the case $n_1=2$,
$\mathfrak{h}_1=\mathfrak{sl}(2,\mathbb{C})$, and
$\mathfrak{h}_2\subsetneq\mathfrak{sp}(n_2,\mathbb{C})$Let
$\mathfrak{k}=\mathfrak{h}_2$. From Proposition \ref{propPkvatK}
below, Lemma \ref{lemTan}, and the considerations of Section
\ref{secRiemannian} it follows that
$\mathfrak{h}=\mathfrak{sp}(1)\oplus\mathfrak{f}\subset\mathfrak{sp}(1)\oplus\mathfrak{sp}(k)\subset\mathfrak{so}(4k)$
is the holonomy algebra of a quaternionic-K\"ahlerian symmetric
space. The theorem is true. $\Box$

\begin{proposition}\label{propPkvatK} Let $\mathfrak{k}\subset\mathfrak{sp}(2m,\mathbb{C})$ be an irreducible
subalgebra, $m\geq 2$. Then
$$\mathcal{P}(\mathfrak{sl}(2,\mathbb{C})\oplus\mathfrak{k})\simeq
\mathbb{C}^2\otimes\mathfrak{g}_1,$$ where $\mathfrak{g}_1$ is the
first Tanaka prolongation of the Lie algebra
$$\mathfrak{g}_{-2}\oplus\mathfrak{g}_{-1}\oplus\mathfrak{g}_0=\mathbb{C} F \oplus\mathbb{C}^{2m}\oplus(\mathfrak{k}\oplus \mathbb{C} H).$$
\end{proposition}

{\bf Proof.} Let $V=\mathbb{C}^{2m}$, let $\Omega$, $\omega$ be
the symplectic forms on $V$ and $\mathbb{C}^2$, and let $e_1$,
$e_2$ be a basis of $\mathbb{C}^2$ such that $\omega(e_1,e_2)=1$.
Let $F,H,E$ be the
 basis of $\mathfrak{sl}(2,\mathbb{C})$ as above. For a linear map $$P:\mathbb{C}^2\otimes
V\to\mathfrak{sl}(2,\mathbb{C})\oplus\mathfrak{k}$$ and $X\in V$
we write
$$P(e_i\otimes X)=\alpha(e_i\otimes X)E+\beta(e_i\otimes
X)F+\gamma(e_i\otimes X)H+T(e_i\otimes X),\quad T(e_i\otimes
X)\in\mathfrak{k},\quad i=1,2.$$ Let us consider the condition
$P\in\mathcal{P}(\mathfrak{sl}(2,\mathbb{C})\oplus\mathfrak{k})$.
Let $X,Y,Z\in V$. Taking the vectors $e_1\otimes X$, $e_1\otimes
Y$, $e_1\otimes Z$, we get
$$\beta(e_1\otimes X)\Omega(Y,Z)+\beta(e_1\otimes Y)\Omega(Z,X)+\beta(e_1\otimes
Z)\Omega(X,Y)=0.$$ Since $\dim V\geq 4$, this implies
$\beta(e_1\otimes X)=0$ for all $X\in V$. Similarly, considering
the vectors $e_2\otimes X$, $e_2\otimes Y$, $e_2\otimes Z$, we get
$\alpha(e_2\otimes X)=0$.

Considering the vectors $e_1\otimes X$, $e_1\otimes Y$,
$e_2\otimes Z$, we obtain
$$\gamma(e_1\otimes X)\Omega(Y,Z)+\Omega(T(e_1\otimes X)Y,Z)-\gamma(e_1\otimes Y)\Omega(X,Z)-\Omega(T(e_1\otimes Y)X,Z)-\beta(e_2\otimes Z)\Omega(Y,X)=0.$$
Let $A\in V$ be the dual vector to $\beta|_{e_2\otimes V}$, i.e.
$\beta(e_2\otimes Z)=\Omega(A,Z)$ for all $Z\in V$. We obtain
$$\gamma(e_1\otimes X)Y+T(e_1\otimes X)Y-\gamma(e_1\otimes Y)X-T(e_1\otimes
Y)X+\Omega(X,Y)A=0.$$ The last equation on $P$ can be obtained in
the same way and it is of the form
$$\gamma(e_2\otimes X)Y-T(e_2\otimes X)Y-\gamma(e_2\otimes Y)X+T(e_2\otimes
Y)X+\Omega(X,Y)B=0,$$ where $B\in V$ is defined by
$\beta(e_1\otimes Z)=\Omega(B,Z)$, $Z\in V$. We conclude that
$P\in\mathcal{P}(\mathfrak{sl}(2,\mathbb{C})\oplus\mathfrak{k})$
if and only if the maps
$$\gamma(e_1\otimes\cdot)H+T(e_1\otimes\cdot),\quad
\gamma(e_2\otimes\cdot)H-T(e_2\otimes\cdot):V\to\mathfrak{k}\oplus\mathbb{C}H
$$ belong to $\mathfrak{g}_1$. Thus,
$$\mathcal{P}(\mathfrak{sl}(2,\mathbb{C})\oplus\mathfrak{k})\simeq\mathfrak{g}_1\oplus\mathfrak{g}_1=\mathbb{C}^2\otimes\mathfrak{g}_1,$$
which is an isomorphism of
$\mathfrak{sl}(2,\mathbb{C})\oplus\mathfrak{k}$-modules. $\Box$

\section{Further remarks}

We are left with the problem to give a direct proof of the fact
that if for a real irreducible representation
$\mathfrak{h}\subset\mathfrak{so}(n)$ of a simple Lie algebra
$\mathfrak{h}$ it holds $\mathcal{P}(\mathfrak{h})\neq 0$, then
$\mathfrak{h}\subset\mathfrak{so}(n)$ is the holonomy algebra of a
Riemannian manifold. The following two cases should be considered:
$\mathcal{P}_0(\mathfrak{h})\neq 0$ and
$\mathcal{P}_1(\mathfrak{h})\neq 0$. It is necessary to prove that
the first condition implies that
$\mathfrak{h}\subset\mathfrak{so}(n)$ is the holonomy algebra of a
not locally symmetric Riemannian manifold; the second condition
implies that $\mathfrak{h}\subset\mathfrak{so}(n)$ may appear as
the holonomy algebra of a symmetric Riemannian manifold. It would
be useful to give a direct proof to the following statement:

{\it If the connected Lie subgroup $H\subset {\rm SO}(n)$
corresponding to an irreducible subalgebra
$\mathfrak{h}\subset\mathfrak{so}(n)$ does not act transitively on
the unite sphere, then $\mathcal{P}_0(\mathfrak{h})=0$.}

The relation of the spaces $\mathcal{P}(\mathfrak{h})$ and
$\mathcal{R}(\mathfrak{h})$ is the following:
$$\mathcal{R}(\mathfrak{h})=\{S\in\mathbb{R}^{n*}\otimes\mathcal{P}(\mathfrak{h})|S(X)(Y)=-S(Y)(X)\}.$$

Consider the natural map
$$\tau:\mathbb{R}^n\otimes\mathcal{R}(\mathfrak{h})\to\mathcal{P}(\mathfrak{h}),\quad
\tau(X\otimes R)=R(X,\cdot)\in\mathcal{P}(\mathfrak{h}).$$ Using
the results of \cite{Leistner}, in \cite{onecomp} it is shown that
$$\tau(\mathbb{R}^n\otimes\mathcal{R}_0(\mathfrak{h}))=\mathcal{P}_0(\mathfrak{h})\,\, (\text{if } n\geq 4),\quad
\tau(\mathbb{R}^n\otimes\mathcal{R}_1(\mathfrak{h}))=\mathcal{P}_1(\mathfrak{h}).$$
It would be useful to get a direct  proof of these statements for
any irreducible subalgebra $\mathfrak{h}\subset\mathfrak{so}(n)$.

Suppose that $\mathcal{P}_1(\mathfrak{h})\neq 0$, i.e.
$\mathcal{P}_1(\mathfrak{h})\simeq \mathbb{R}^n$. Then there
exists an $\mathfrak{h}$-equivariant linear isomorphism
$S:\mathbb{R}^n\to \mathcal{P}_1(\mathfrak{h})$ defined up to a
constant multiple. It should be proved that $S(X)(Y)=-S(Y)(X)$,
i.e. $S\in\mathcal{R}_1(\mathfrak{h})$.

 The space $\mathcal{P}(\mathfrak{h})$ is contained in the
tensor product $\mathbb{R}^n\otimes\mathfrak{h}$.  A statement
form \cite{onecomp} implies that the decomposition of
$\mathbb{R}^n\otimes\mathfrak{h}$ into the sum of irreducible
$\mathfrak{h}$-modules is of the form
$$\mathbb{R}^n\otimes\mathfrak{h}= k\mathbb{R}^n\oplus (\oplus_\lambda V_\lambda),$$ where
$k$ is the number of non-zero labels on the Dynkin diagram for the
representation of $\mathfrak{h}\otimes\mathbb{C}$ on
$\mathbb{C}^n$, and $V_\lambda$ are pairwise non-isomorphic
irreducible $\mathfrak{h}$-modules that are not isomorphic to
$\mathbb{R}^n$. If $\mathcal{P}_0(\mathfrak{h})\neq 0$, then it
coincides with the highest irreducible component in
$\mathbb{R}^n\otimes\mathfrak{h}$.

The space $\mathcal{R}(\mathfrak{h})$ is contained in
$\odot^2\mathfrak{h}$ \cite{Al}. If
$\mathcal{R}_1(\mathfrak{h})\neq 0$, then it is spanned by the map
$\id_{\mathfrak{h}}\in\odot^2\mathfrak{h}\subset\wedge^2\mathbb{R}^n\otimes\mathfrak{h}$,
note that $\id_{\mathfrak{h}}(X,Y)=\pr_\mathfrak{h}(X\wedge Y)$.
Consequently, if $\mathcal{R}_1(\mathfrak{h})\neq 0$, then
$$\mathcal{P}_1(\mathfrak{h})=\tau(\mathbb{R}^n\otimes\id_\mathfrak{h})=\{\pr_\mathfrak{h}(X\wedge\cdot)|X\in\mathbb{R}^n\}.$$
But it is not clear why if $\mathcal{P}_1(\mathfrak{h})\neq 0$,
then it should coincide with
$\tau(\mathbb{R}^n\otimes\id_\mathfrak{h})$ (such statement would
imply $\mathcal{R}_1(\mathfrak{h})\simeq\mathbb{R}$).

The statement of the following lemma can be checked directly.

\begin{lemma} Let $S:\mathbb{R}^n\to\mathcal{P}(\mathfrak{h})$ be a linear map. Consider the
map $$T:\wedge^2\mathbb{R}^n\to\mathfrak{h},\quad
T(X,Y)=S(X)(Y)-S(Y)(X).$$ Then
$T+T^*\in\mathcal{R}(\mathfrak{so}(n))$, where
$T^*:\mathfrak{so}(n)\to\mathfrak{so}(n)$ is given by
$(T^*(X,Y)Z,W)=(T(Z,W)X,Y)$.\end{lemma}

We are able to show that the condition
$\mathcal{P}_1(\mathfrak{h})\neq 0$ implies
$\mathcal{R}_1(\mathfrak{h})\neq 0$ only under an additional
assumption on the representation
$\mathfrak{h}\subset\mathfrak{so}(n)$.

\begin{proposition} Let $\mathfrak{h}\subset\mathfrak{so}(n)$ be an irreducible
representation of real type of a simple Lie algebra $\mathfrak{h}$
such that $\mathcal{P}_1(\mathfrak{h})\neq 0$. If the irreducible
representation
$\mathfrak{h}\otimes\mathbb{C}\subset\mathfrak{so}(n,\mathbb{C})$
is given by the Dynkin diagram with only one or two non-zero
labels, then $\mathcal{R}_1(\mathfrak{h})\neq 0$, i.e.
$\mathfrak{h}\subset\mathfrak{so}(n)$ is the holonomy algebra of a
symmetric Riemannian space.
\end{proposition}

{\bf Proof.} If the non-zero label is only one, then the
multiplicity of $\mathbb{R}^n$ in the tensor product
$\mathbb{R}^n\otimes\mathfrak{h}$ is one, namely the submodule of
$\mathbb{R}^n\otimes\mathfrak{h}$ isomorphic to $\mathbb{R}^n$ is
equal to $\tau(\mathbb{R}^n\otimes\id_\mathfrak{h})$, this implies
the proof.

Suppose that there are two non-zero labels. Then the multiplicity
of $\mathbb{R}^n$ in the tensor product
$\mathbb{R}^n\otimes\mathfrak{h}$ is two. One submodule isomorphic
to $\mathbb{R}^n$ is equal to
$\tau(\mathbb{R}^n\otimes\id_\mathfrak{h})$. The orthogonal
complement to $\tau(\mathbb{R}^n\otimes\id_\mathfrak{h})$ in
$\mathbb{R}^n\otimes\mathfrak{h}$  is the subspace
$(\mathbb{R}^n\otimes\mathfrak{h})_0\subset
\mathbb{R}^n\otimes\mathfrak{h}$ consisting of linear maps
$\varphi:\mathbb{R}^n\to\mathfrak{h}$ with $\tRic(\varphi)=0$
\cite{onecomp}. This space contains a uniquely defined submodule
isomorphic to $\mathbb{R}^n$. It is obvious that the projection of
$\mathcal{P}_1(\mathfrak{h})$ to
$\tau(\mathbb{R}^n\otimes\id_\mathfrak{h})$ is not trivial.
Clearly, the subspace
$W\subset\mathbb{R}^n\otimes\mathbb{R}^n\otimes\mathfrak{h}$ of
elements annihilated by $\mathfrak{h}$ is two-dimensional; it
contains the subspace
$\mathbb{R}\id_\mathfrak{h}\subset\odot^2\mathfrak{h}\subset\wedge^2\mathbb{R}^n\otimes\mathfrak{h}$.
Since $\mathcal{P}_1(\mathfrak{h})\simeq \mathbb{R}^n$, there
exists an $\mathfrak{h}$-equivariant isomorphism
$S:\mathbb{R}^n\to\mathcal{P}_1(\mathfrak{h})$, $S\in W$. If
$W\subset \wedge^2\mathbb{R}^n\otimes\mathfrak{h}$, then
$S\in\mathcal{R}_1(\mathfrak{h})$. Otherwise,
$W=\mathbb{R}\id_\mathfrak{h}\oplus\mathbb{R}\psi$, where
$\psi\in\odot^2\mathbb{R}^n\otimes\mathfrak{h}$. Since
$\mathcal{P}_1(\mathfrak{h})\not\subset(\mathbb{R}^n\otimes\mathfrak{h})_0$,
$S\notin\mathbb{R}\psi$. The element
$T\in\wedge^2\mathbb{R}^n\otimes\mathfrak{h}$ defined in the above
lemma belongs to $W$, hence $T=c\id_\mathfrak{h}$ for some
non-zero $c\in\mathbb{R}$. Next,
$\id_\mathfrak{h}^*=\id_\mathfrak{h}$, and from the lemma it
follows that $\id_\mathfrak{h}\in\mathcal{R}_1(\mathfrak{h})$,
i.e. $\mathcal{R}_1(\mathfrak{h})=\mathbb{R}\id_\mathfrak{h}$.
This proves the proposition. $\Box$

%\section*{References}

\end{document}